\documentclass[12pt]{article}

\usepackage{amsmath}
\usepackage{amsfonts}
\usepackage{amssymb}
\usepackage{amsthm}
\usepackage{diagrams}

\newtheorem{lem}{Lemma}[section]
\newtheorem{cor}[lem]{Corollary}
\newtheorem{prop}[lem]{Proposition}
\newtheorem{thm}[lem]{Theorem}
\newtheorem{conj}[lem]{Conjecture}
\newtheorem{Defn}[lem]{Definition}
\newtheorem{Ex}[lem]{Example}
\newtheorem{Question}[lem]{Question}

\newenvironment{ex}{\begin{Ex}\rm}{\end{Ex}}
\newenvironment{question}{\begin{Question}\rm}{\end{Question}}

\newcommand{\ideal}[1]{\mathfrak{#1}}

\newcommand{\m}{\ideal{m}}
\newcommand{\n}{\ideal{n}}
\newcommand{\p}{\ideal{p}}
\newcommand{\q}{\ideal{q}}

\newcommand{\s}{\ideal{s}}

\newcommand{\red}[1]{{#1}_{\text{red}}}

\newcommand{\A}{\mathbb{A}}

\newcommand{\ol}{\overline}

\newcommand{\x}{\mathbf{x}}

\newcommand{\rank}{\mathrm{rank}}

\newcommand{\Ht}{\mathrm{ht}\,}	

\newcommand{\spec}{\text{Spec}}

\newcommand{\ass}{\text{Ass}}	
\newcommand{\len}{\ell}
\newcommand{\tor}{\text{Tor}}

\newcommand{\gr}{\mathrm{gr}}

\newarrowtail{<=}\Leftarrow\Rightarrow{\@cmex7F}{\@cmex7E}

\newarrow{Implies}{=}{=}{=}{=}{=>}
\newarrow{Iff}{<=}{=}{=}{=}{=>}

\begin{document}
\bibliographystyle{amsplain}
\begin{center}
\Large INTERSECTIONS OF SYMBOLIC POWERS OF PRIME IDEALS

\medskip
\large SEAN SATHER-WAGSTAFF\footnotetext{2000 Mathematics Subject 
Classification 13H05, 13H15, 13C15, 13D22}

\bigskip
\normalsize\textsc{Abstract}
\end{center}

\footnotesize Let $(R,\m)$ be a local ring with prime ideals $\p$ and $\q$ such 
that $\sqrt{\p+\q}=\m$.  If $R$ is regular and contains a field, and 
$\dim(R/\p)+\dim(R/\q)=\dim(R)$, we prove that 
$\p^{(m)}\cap\q^{(n)}\subseteq\m^{m+n}$ for all positive 
integers $m$ and $n$.  
This is proved using a generalization of Serre's 
Intersection Theorem which we apply to a hypersurface $R/fR$.  The generalization 
gives conditions that guarantee that Serre's bound on the intersection 
dimension $\dim(R/\p)+\dim(R/\q)\leq\dim(R)$ holds when $R$ is nonregular.  


\normalsize 
\section{Introduction} \label{sec:intro}
Let $k$ be an algebraically closed 
field of characteristic 0 and $\A^d_k$ the $d$-dimensional affine space 
over $k$.  That is, $\A^d_k=\spec(S)$ where $S$ is the polynomial 
algebra $k[X_1,\ldots,X_d]$.  
Let $Y$ and $Z$ be closed subvarieties of $\A^d_k$ that 
intersect at a finite number of points, and assume without loss of 
generality that the intersection contains the origin.  Let $f$ be a nonzero 
regular function on $\A^n_k$ that vanishes on $Y\cup Z$, and let $m$ and $n$ 
denote the orders of vanishing of $f$ along $Y$ and $Z$, respectively.  
That is, all the partial derivatives of $f$ of order at most
$(m-1)$ vanish 
identically on $Y$ and at least one $m$th order partial does not, 
and similarly for $n$.  
It is reasonable to ask, in this case, what 
the order of vanishing of $f$ is at the origin.  
Of course, if 
$r=\max\{m,n\}$ then $f$ vanishes to order at least $r$, and
the following easy example shows that, without additional hypotheses, 
this bound is sharp.

\begin{ex} \label{ex:trivial}
Let $d=3$, $Y=V(X_1,X_2)$ and $Z=V(X_2,X_3)$.  Then $Y$ and $Z$ intersect 
precisely at the origin, and the function $X_2$ 
vanishes exactly to order 1 at every point of $Y$ and at every point of 
$Z$.  In particular, $X_2$ vanishes exactly to order 1 at the origin.
\end{ex}

One may feel somewhat cheated by this example since $X_2$ is an 
``essential'' defining polynomial for both subvarieties.  
The existence of this polynomial is objectionable for another related reason.  In 
our example, we are considering two lines in 3-space, and in general, two such 
varieties will not intersect.  In fact, our lines intersect exactly because 
they are coplanar, and more specifically because they live in the plane defined by the 
vanishing of the ``offending'' polynomial $X_2$.  If we consider two lines in general 
position, that is if we
``wiggle'' the lines slightly, then they will generally be pulled
away from one 
another;  the question of vanishing of functions at an intersection 
point is then meaningless.  This suggests that we should demand that the
objects under consideration be 
well-behaved under ``wiggling.''  In order to guarantee this, we 
recall a classical theorem from algebraic geometry which states that, under our 
initial hypotheses, $\dim(Y)+\dim(Z)\leq d$.  We shall refer to this 
as the Classical Intersection Theorem.
In addition, when $\dim(Y)+\dim(Z)=d$, the varieties will intersect 
generically at a finite nonempty set of points.
In Example~\ref{ex:trivial} this 
inequality is strict, and we can directly attribute this to the smooth 
hypersurface 
$V(X_2)$ containing $Y$ and $Z$.  
What can be said about orders of vanishing 
if we assume that the subvarieties satisfy the additional hypothesis 
$\dim(Y)+\dim(Z)=d$?
To motivate our answer, we consider the following example.

\begin{ex} \label{ex:transverse}
Assume that $d$ is at least 2 and fix an integer $i$ between 1 and $d-1$.  
Let $\p=(X_1,\ldots,X_i)S$, $\q=(X_{i+1},\ldots, X_d)S$, $Y=V(\p)$ 
and $Z=V(\q)$.  Then $Y$ and $Z$ intersect at exactly the origin and 
$\dim(Y)+\dim(Z)=d$.  Zariski's Main Lemma on Holomorphic 
Functions~\cite{zariski:flthfav} 
implies that a nonzero 
regular function $f$ vanishes to order $m$ at 
every point of $Y$ if and only if $f\in\p^{(m)}$, where 
$\p^{(m)}=R\cap \p^m R_{\p}$ is the 
$m$th symbolic power of $\p$.  
In this case, the ordinary and symbolic 
powers of $\p$ agree, so $f$ vanishes to order $m$ at 
every point of $Y$ if and only if $f\in\p^m$.  Similarly, $f$ vanishes to 
order $n$ at every point of $Z$ if and only if $f\in\q^{(n)}=\q^n$.  
It is straightforward to show that, in this case
\[ \p^{(m)}\cap\q^{(n)}=\p^m\cap\q^n=\p^m\q^n\subseteq\m^{m+n}. \]
Therefore $f$ vanishes at the origin 
to order at least $m+n$.  That is, the sum of the 
orders of vanishing along $Y$ and $Z$ gives a lower bound for the order of 
vanishing at the origin.  One can easily construct examples to show that 
this bound is sharp.
\end{ex}

This example leads us to ask the following question.

\begin{question} \label{question:affine}
Let $k$ be an algebraically closed field of characteristic 0, and $Y$ and $Z$  
closed subvarieties of $\A^d_k$ that intersect at 
finitely many points, including the origin, and such that 
$\dim(Y)+\dim(Z)=d$.  If $f$ is a nonzero regular 
function on $\A^n_k$ vanishing along $Y$ and $Z$ to orders $m$ and $n$, 
respectively, does $f$ vanish to order at least $m+n$ at the origin?
\end{question}

As a corollary to one of the main results of this paper, we answer this 
question in the affirmative.  See Corollary~\ref{cor:affine} below.  

The discussion preceding in Example~\ref{ex:transverse} suggests a purely
algebraic formulation of Question~\ref{question:affine}.  We state the 
local version as part of Conjecture~\ref{conj:all} below.  Before we do so, 
we present a somewhat more algebraic motivation for 
this question.  We include it, not only
because the methods we use to 
answer Question~\ref{question:affine} are purely algebraic in nature,
but also because it indicates how we originally 
came to consider these ideas.

Let $(R,\m)$ be a local, Noetherian ring with prime ideals $\p$ 
and $\q$ such that $\sqrt{\p+\q}=\m$.   
Serre~\cite{serre:alm} generalized the Classical Intersection Theorem 
by proving that, if $R$ is regular, then 
\[ \dim(R/\p)+\dim(R/\q)\leq\dim(R).\]
We shall refer to this result as Serre's Intersection Theorem.    
When $R/\p$ has finite projective dimension (true automatically 
if $R$ is regular) Serre defined the intersection multiplicity of $R/\p$ 
and $R/\q$ as
\[ \chi(R/\p,R/\q)=\!\!\sum_{i=0}^{\dim(R)} \!\!(-1)^i\len(\tor^R_i(R/\p,R/\q))\]
where $\ell$ is the length function.
When $R$ is regular and unramified, he proved that $\chi(R/\p,R/\q)\geq 0$ 
with strict inequality holding if and only if 
$\dim(R/\p)+\dim(R/\q)=\dim(R)$.  He conjectured that the same holds true 
if the ring is ramified.

Independently, Roberts~\cite{roberts:vimpc} and 
Gillet-Soul{\'e}~\cite{gillet:knmi} proved the vanishing part of the conjecture in the 
ramified case:  if 
$\dim(R/\p)+\dim(R/\q)<\dim(R)$ then $\chi(R/\p,R/\q)=0$.  
Gabber~\cite{berthelot:ava, hochster:nimrrlr, roberts:rdsmcgpnc} 
proved nonnegativity:  $\chi(R/\p,R/\q)\geq 0$.  The Positivity 
Conjecture is the converse to the vanishing result and is still 
open in the ramified case.  

By applying Gabber's methods to the Positivity Conjecture, 
Kurano and Roberts~\cite{kurano:pimsppi} proved the following.

\begin{thm} \label{thm:KR}
Let $(R,\m)$ be a regular local ring with prime ideals $\p$ and 
$\q$ such that $\sqrt{\p+\q}=\m$ and $\dim(R/\p)+\dim(R/\q)=\dim(R)$.  If
$R$ is ramified or contains a field and $\chi(R/\p,R/\q)>0$, 
then $\p^{(m)}\cap\q\subseteq\m^{m+1}$ for all positive integers $m$, 
where $\p^{(m)}$ is the $m$th symbolic power of $\p$.
\end{thm}

Because Kurano and Roberts expect the Positivity Conjecture to be true, 
this result motivated them to make the following conjecture.

\begin{conj} \label{conj:KR}
Let $(R,\m)$ be a regular local ring with prime ideals $\p$ and 
$\q$ such that $\sqrt{\p+\q}=\m$ and $\dim(R/\p)+\dim(R/\q)=\dim(R)$.  
Then $\p^{(m)}\cap\q\subseteq\m^{m+1}$ for all positive integers $m$.
\end{conj}

There is a significant amount of interest in 
the properties of symbolic powers of prime ideals in regular local rings.
In addition to the paper of Kurano-Roberts, the interested reader should 
refer to~\cite{eisenbud:essfi, hochster:csopi, huckaba:pihsad, 
huneke:ssrlr}.

Theorem~\ref{thm:KR} shows that Conjecture~\ref{conj:KR} is true for 
regular rings containing an arbitrary field since positivity holds for 
these rings.  
In their paper, Kurano and Roberts verify the conjecture, with no reference 
to positivity, 
when $R$ contains a field of 
characteristic 0.  They prove this by first demonstrating that, if   
$\p\subseteq\m^2$ then $\p^{(m)}\subseteq\m^{m+1}$ for all $m\geq 1$;
they then reduce the general case to the case where $\p\subseteq\m^2$.  
Curiously, the proof makes little reference to the second prime $\q$.
The current author~\cite{sather:dicmr} has verified the conjecture, again 
with no reference to positivity,  when 
$R$ contains a field of arbitrary characteristic.  This is accomplished 
by first proving a 
generalization of Serre's Intersection Theorem and applying it to a 
hypersurface $R/fR$.  We 
discuss the interaction between such generalizations and 
Conjecture~\ref{conj:KR} (and its generalizations) in greater depth 
below.  We note that Conjecture~\ref{conj:KR} is still open 
in mixed-characteristic.   

Theorem~\ref{thm:KR} is immediately applicable to 
Question~\ref{question:affine}, as follows.  Let $Y$ and $Z$ be closed subvarieties 
of $\A_k^d$ that intersect in a finite number of points, including the 
origin, and such that $\dim(Y)+\dim(Z)=d$.  Let $f$ be a nonzero function 
on $\A^n_k$ that vanishes to order $m$ and $n$ along $Y$ and $Z$, 
respectively.  Without loss of 
generality, assume that $m\geq n$.  Let
$\p$ and $\q$ be the prime ideals of $S$ defining $Y$ and $Z$, 
respectively, and let $\m=(X_1,\ldots,X_d)S$.  
The function $f$ is an element of $\p^{(m)}$ by Zariski's 
Main Lemma, and it is an element of $\q$ by assumption.  It follows from
Theorem~\ref{thm:KR} that 
\[ f\in S\cap\m^{m+1}S_{\m}=\m^{(m+1)}=\m^{m+1}.\]  
In other words, $f$ 
vanishes to order at least $m+1$ at the origin.  Although this does not 
answer Question~\ref{question:affine}, it at least demonstrates that the assumption 
$\dim(Y)+\dim(Z)=d$ implies greater order of vanishing than demonstrated in 
Example~\ref{ex:trivial}.  

In considering Conjecture~\ref{conj:KR} with
Question~\ref{question:affine} in mind, one might be tempted to ask why the 
conjecture is not more symmetric.  (This is, of course, an easy question 
to ask in light of 
the discussion above.  However, one might also be led to this question by 
observing the symmetry of the intersection multiplicity $\chi$ and its 
relevance to the conjecture via Theorem~\ref{thm:KR}.)
The main reason is that the evidence that led 
to the formulation of 
the conjecture was Theorem~\ref{thm:KR}, and at the time there was no 
similar result suggesting stronger results for
$\p^{(m)}\cap\q^{(n)}$.   In addition, one notices that the proof of 
Theorem~\ref{thm:KR} is not symmetric, as one of the first steps is to 
take a regular alteration of $\spec(R/\p)$\footnote{A regular alteration of 
$\spec(R/\p)$ is similar to a 
resolution of singularities:  it is a projective morphism 
$\phi:X\rightarrow\spec(R/\p)$ where $X$ is a regular scheme.  Regular 
alterations are weaker than 
resolutions of singularities because they do not necessarily induce 
isomorphisms on the fields of rational functions, only finite extensions.  
However, unlike resolutions of singularities, they are known to exist for 
rings essentially of finite type over a field (of arbitrary characteristic) 
or a complete discrete valuation ring, by the work of
de Jong~\cite{dejong:ssa}.  The existence of regular alterations is key to 
Gabber's proof of the nonnegativity 
conjecture.}
and pass to a closely related projective scheme.
Therefore, Kurano and Roberts did not necessarily expect a 
symmetric result.

As the earlier discussion indicates, though, we do consider a more 
symmetric intersection in the current paper.
More specifically,  under the hypotheses of 
Conjecture~\ref{conj:KR}, we ask whether the containment 
$\p^{(m)}\cap\q^{(n)}\subseteq\m^{m+n}$ holds.
In one of our main results 
we answer this question in the affirmative for rings 
containing an arbitrary field.  
More specifically, we have the following.

\bigskip
\noindent\textbf{Theorem~\ref{thm:powers}}
\textit{Let $(R,\m)$ be a regular local ring containing a field and 
$\p$ and $\q$ prime ideals of $R$ such that $\sqrt{\p+\q}=\m$ and 
$\dim(R/\p)+\dim(R/\q)=\dim(R)$.  Then 
$\p^{(m)}\cap\q^{(n)}\subseteq\m^{m+n}$ for all $m,n\geq 1$.}

\bigskip

As an immediate consequence of this result, we are able to answer 
Question~\ref{question:affine} in the affirmative.

Theorem~\ref{thm:powers} is established by proving the following
generalization of Serre's Intersection Theorem for nonregular 
rings and applying it to the hypersurface $R/fR$ where $f$ is an 
element of $\p^{(m)}\cap\q^{(n)}$.  

\bigskip
\noindent\textbf{Theorem~\ref{thm:id-1}}
\textit{Let $(A,\n)$ be a quasi-unmixed local ring
and assume that one of the following conditions holds.
\begin{enumerate}
\item $\red{A}$ contains a field, or
\item $A$ is a ring of mixed-characteristic such that the residual 
characteristic $p$ is 
part of a reductive system of parameters of $A$.
\end{enumerate}
Let $P$ and $Q$ be prime ideals in $A$ such 
that both $A/P$ and $A/Q$ are analytically unramified,
$\sqrt{P+Q}=\n$, and $e(A)<e(A_P)+e(A_Q)$.  Then 
$\dim(A/P)+\dim(A/Q)\leq\dim(A)$.}

\bigskip

When we say that $A$ is quasi-unmixed, we mean that every irreducible component of the 
completion $A^*$ has the same dimension, that is, that $A^*$ is 
equidimensional.  By a theorem of Ratliff~\cite{ratliff:qldafccpiII}, this is 
equivalent to $A$ being equidimensional and universally catenary.
The ring $\red{A}$ is the reduced ring $A/\text{nil}(A)$, 
and $e(A)$ is the Hilbert-Samuel multiplicity of $A$.  
A reductive system of parameters for $A$ is a system of parameters that 
generates a reduction of $\n$.  (See Section~\ref{sec:equalcharacteristic} for 
complete definitions and relevant
properties.)  
Briefly, the 
connection between Theorems~\ref{thm:id-1} and~\ref{thm:powers} is the fact 
that, if $0\neq f\in\p$, then $e(R_{\p}/fR_{\p})=m$ if and only if 
$f\in\p^{(m)}\smallsetminus\p^{(m+1)}$.  

This leads us to consider two conjectural generalizations of Serre's 
Intersection Theorem.  Each statement generalizes a 
question about symbolic powers of prime ideals in regular local rings.  We 
list the conjectures here with the corresponding conjectures for symbolic 
powers.  

\begin{conj} \label{conj:all}
Let $(A,\n)$ be a quasi-unmixed local ring with prime ideals $P$ and 
$Q$ such that $\sqrt{P+Q}=\n$.  Also, let $(R,\m)$ be a regular local 
ring with prime ideals $\p$ and $\q$ such that $\sqrt{\p+\q}=\m$ and 
$\dim(R/\p)+\dim(R/\q)=\dim(R)$.  Then
\begin{description}
\item[(SP-1)]  $\p^{(m)}\cap\q\subseteq\m^{m+1}$ for all $m\geq 1$.
\item[(SP-2)]  $\p^{(m)}\cap\q^{(n)}\subseteq\m^{m+n}$ for all $m,n\geq 1$.
\item[(ID-1)]  If $e(A)=e(A_P)$ and $A/P$ is analytically unramified, then 
\[ \dim(A/P)+\dim(A/Q)\leq\dim(A).\]
\item[(ID-2)]  If $e(A)<e(A_{P})+e(A_{Q})$ and both $A/P$ and $A/Q$ are 
analytically unramified,  then 
\[ \dim(A/P)+\dim(A/Q)\leq\dim(A).\]
\end{description}
\end{conj}

We immediately observe that (ID-1) and (ID-2) are conjectural
generalization of Serre's Intersection Theorem.  This is due to the fact 
that, if $A$ is regular then $e(A)=e(A_P)=e(A_Q)=1$.  
Conjecture (SP-1) and (ID-1) have 
been verified for rings containing an arbitrary field and 
for a number of special cases in~\cite{kurano:pimsppi, sather:dicmr, 
sather:mdium}, and 
Theorems~\ref{thm:powers} and~\ref{thm:id-1} of the present paper
imply that (SP-2) and 
(ID-2) are true for rings containing fields. 

We note that ours are not the first efforts made to generalize Serre's 
Intersection Theorem to nonregular rings.  Of course, there is the landmark 
paper of Peskine and Szpiro~\cite{peskine:sm}, as well as the more recent 
work of Simon~\cite{simon:idcm}.  

Below, we discuss the reasons 
for the technical assumptions in Conjecture~\ref{conj:all}.  Before we do so, 
we discuss the relations between the individual conjectures.
Obviously, (SP-2) implies (SP-1).
If (ID-2) holds for all complete hypersurfaces, then (SP-2) 
holds for all regular local rings.  To see this, pass to the completion 
$R^*$ of $R$ and apply (ID-2) as in the proof of Theorem~\ref{thm:powers}.  
Similarly, if (ID-1) is true for all complete hypersurfaces, then (SP-1) 
is true for all regular local rings. 
Clearly (ID-2) implies (ID-1) when $A/Q$ is analytically unramified, 
because $e(A_Q)>0$ and 
$e(A_P)=e(A)$ imply that $e(A)<e(A_P)+e(A_Q)$.  
It is straightforward to reduce (ID-1) 
to the case of a complete domain---we perform this reduction explicitly 
for (ID-2) in the proof of Theorem~\ref{thm:id-1}---and it follows that (ID-2) 
implies (ID-1). 
We summarize the implications 
is the following diagram.
\begin{diagram}[small]
\text{(ID-2)} & \rImplies & \text{(SP-2)} \\
\dImplies & & \dImplies \\
\text{(ID-1)} & \rImplies & \text{(SP-1)}
\end{diagram}

By Theorem~\ref{thm:KR}, Serre's Positivity Conjecture implies (SP-1) for 
ramified rings.  At this time, we do not know whether positivity implies 
(SP-1) for unramified rings or any of the other conjectures in any case of 
mixed characteristic.  In addition
we do not know whether any of these conjectures 
imply Serre's Positivity Conjecture.

Here we discuss the technical assumptions of Conjecture~\ref{conj:all}.
It is easy to find examples to show why $R$ must be regular, 
$\sqrt{\p+\q}$ must equal $\m$, and $\dim(R/\p)+\dim(R/\q)$ must equal 
$\dim(R)$.  Similarly, one sees that $\sqrt{P+Q}$ must equal $\n$.  
The ring $A$ must be equidimensional in (ID-1) and (ID-2), by 
considering  $A=k[\![X,Y,Z]\!]/(XY,XZ)$ and letting $P$ and 
$Q$ be the minimal primes of $A$.  In~\cite{sather:mdium} Example 6.5, we give
an example showing that $A$ should be at least 
catenary.  
We may or may not need the full strength of quasi-unmixed, that is, 
equidimensional and universally catenary.
In our arguments, where we pass to the completion, we need the 
completion to be equidimensional.
If $A$ is excellent and equidimensional, 
then it is automatic, so we do not consider this 
assumption too restrictive.  

When we wish to use $e(A_P)$, we assume that $A/P$ 
is analytically unramified.  By this, we mean 
that the completion $(A/P)^*=A^*/PA^*$ is reduced, or in other words, the 
ideal $PA^*$ 
is an intersection of prime ideals of $A^*$.
The purpose of this assumption is to guarantee that 
the multiplicity is well-behaved under localization.  More specifically, we 
require that the multiplicity not increase after localizing.  Our 
guarantee is from Nagata~\cite{nagata:lr} (40.1).

\begin{thm} \label{thm:nagata}
Let $P$ be a prime ideal of a local ring $A$.  If $A/P$ is analytically 
unramified and if $\Ht(P)+\dim(A/P)=\dim(A)$, then $e(A_P)\leq e(A)$.
\end{thm}

Regarding the analytically unramified assumption,
Nagata~\cite{nagata:lr} (Appendix A2) wrote the following. 
``It is not yet known to the writer's knowledge 
whether or not (40.1) is true without assuming that $P$ is analytically 
unramified.''  (The statement ``$P$ is analytically unramified'' 
means ``$A/P$ is analytically unramified.'')
If it is shown that this condition can be omitted from the 
statement of Theorem~\ref{thm:nagata}, then the corresponding conditions 
should probably be omitted from Conjectures (ID-1) and (ID-2).  For 
now, however, we leave them intact, especially because our arguments (in 
which we pass to the completion) depend on the assumptions.  As with the 
quasi-unmixedness assumption, if $A$ is excellent then $A/P$ is guaranteed 
to be analytically unramified, so this is not an unbearable restriction.


Acknowledgments:  I am grateful to S.~Dutta, P.~Griffith, G.~Leuschke 
and A.~Singh for their helpful comments and suggestions.

\section{The Main Results} \label{sec:equalcharacteristic}

Before we verify the conjectures for rings containing fields, we give some 
definitions and background results.  An excellent reference for 
multiplicity-theoretic results is Herrmann, Ikeda and 
Orbanz~\cite{herrmann:ebu}.

If $(A,\n)$ is a local ring and $M$ is a 
nonzero, finitely generated $A$-module of dimension $d$ 
and $I$ is an ideal of $A$ such that $M/IM$ 
has finite length, then the Hilbert function $P(n)=\len_A(M/I^{n+1}M)$ is 
a polynomial of degree $d$ for $n\gg 0$.  The Hilbert-Samuel multiplicity 
is the positive integer $e_A(I,M)$ such that
$P(n)=\frac{1}{d!}e_A(I,M)n^d+(\text{lower degree terms})$.  When there is no 
danger of confusion, we shall write $e(I,M)$.  If $I=\n$, we write 
$e_A(M)$ or $e(M)$.

Regarding the Hilbert-Samuel multiplicity,  
we shall use the following facts freely and possibly without reference.

\begin{enumerate}
\item (Associativity Formula)
Let $A$ be a local ring with finitely generated module $M$.  Then
\begin{equation} 
e_A(M)=\sum_{P}\len(M_P)e(A/P) \label{eqn:associativity}
\end{equation}
where the sum is taken over all prime ideals $P$ in the support of $M$ such that $\dim(A/P)=\dim(M)$.  
In particular, this sum is finite.

\item Let $(R,\m)$ be a local ring contained in a ring 
$A$ such that the extension $R\hookrightarrow A$ is module-finite.
Then $A$ is a semilocal ring, say with maximal ideals $\n_1,\ldots,\n_n$, 
such that $\dim(A)=\dim(R)$ and 
each $\n_i\cap R=\m$.   Then
\begin{equation} 
e_{R}(\m,A)=
  \sum_i [A/\n_i:R/\m] e_{A_{\n_i}}\!(\m A_{\n_i},A_{\n_i}) 
    \label{eqn:mult}
\end{equation}
where the sum is taken over all indices $i$ such that $\Ht(\n_i)=\dim(A)$.  

\item Let $R$ be a Noetherian ring contained in a ring $A$
such that the extension $R\hookrightarrow A$ is module-finite, and 
let $\s$ be a prime ideal 
of $R$.  Let $\{ S_1,\ldots,S_j\}$ be the set of prime ideals of $A$ such that 
$S_i\cap R=\s$ and $\Ht(S_i)=\Ht(\s)$.  Then
\begin{equation}
e_{R_{\s}}(\s R_{\s},A_{\s})=
  \sum_{i=1}^j [\kappa(S_i):\kappa(\s)] e_{A_{S_i}}\!(\s A_{S_i},A_{S_i})
     \label{eqn:mult2}
\end{equation}
where $\kappa(\s)$ is the residue field of $R_{\s}$ and similarly for 
$\kappa(S_i)$.

\item Let $A\rightarrow A'$ be a flat, local 
homomorphism of local rings $(A,\n)$ 
and $(A',\n')$ such that $\n A'=\n'$.  Then $e(A')=e(A)$.

\item Let $A$ be a local domain and 
$I$ an ideal of $A$ such that $A/I$ has finite length.  If $M$ is an 
$A$-module of positive rank $r$, then 
\begin{equation}
e_A(I,M)=e_A(I,A)\cdot r. \label{eqn:rank}
\end{equation}
\end{enumerate}

Property 1 is Bruns and Herzog~\cite{bruns:cmr} Corollary 4.7.8.
Property 2 is Nagata~\cite{nagata:lr} (23.1).  Property 3 is simply 
(\ref{eqn:mult}) applied to the extension $R_{\s}\rightarrow A_{\s}$.
Property 4 is Herzog~\cite{herzog:omlr} Lemma 2.3.
And property 5 is Matsumura~\cite{matsumura:crt} Theorem 14.8.

The following is the most important tool used to prove out main results.

\begin{lem} \label{lem:newrank}
Let $(A,\n)$ be an equidimensional, catenary local ring containing a
catenary local domain $(R,\m)$ such that the extension 
$R\hookrightarrow A$ is module-finite.  Let $P$ and $Q$ be prime ideals of $A$ 
such that $\sqrt{P+Q}=\n$ and both $A/P$ and $A/Q$ are analytically 
unramified,
and let $\p=P\cap R$ and $\q=Q\cap R$.  Assume that the quotient $R/\s$ is 
analytically unramified for every prime 
ideal $\s\supset\p+\q$. If 
$e_R(\m,A)<e(A_P)+e(A_Q)$, then $\sqrt{\p+\q}=\m$.  
\end{lem}

If $R$ is an excellent local domain and $A$ is equidimensional, then $A$ is 
also excellent.  In this case, $A$ is automatically catenary and the 
quotients $A/P$, $A/Q$ and $R/\s$ are guaranteed to be analytically unramified.
In our applications, both $A$ and $R$ will in fact be complete, so 
the reader is free to assume either of these stronger conditions.

\begin{proof}
Let $\s$ be a prime ideal of $R$ such that $\p+\q\subseteq\s$.  To show 
that $\sqrt{\p+\q}=\m$, it suffices to show that $\s=\m$.  By the going-up 
property for integral extensions, fix prime ideals $S_1$ and $S_2$ of $A$ 
such that $S_1\supseteq P$, $S_2\supseteq Q$ and $S_1\cap R=S_2\cap R=\s$.  
To show that $\s=\m$, it suffices to show that $S_1=S_2$ because then 
$S_1\supseteq P+Q$ so that $S_1=\n$ and $\s=S_1\cap R=\n\cap R=\m$.

Suppose that $S_1\neq S_2$.  We first observe that, for every prime $S$ of 
$A$ such that $S\cap R=\s$, we have $\Ht(S)=\Ht(\s)$.  This is due to the 
fact that, by integrality $\dim(R/\s)=\dim(A/S)$, and the
catenary and equidimensional assumptions imply that $\Ht(S)=\dim(A)-\dim(A/S)$ and 
similarly for $\s$.
In particular, $\Ht(S_i)=\Ht(\s)$ for $i=1,2$.  Thus
\begin{align*}
e(A_P)+e(A_Q)
  & > e_R(\m,A)=\rank_R(A)e(R) & \text{(by (\ref{eqn:rank}))} \\
  & \geq \rank_{R_{\s}}(A_{\s})e(R_{\s}) & \text{(by 
         Theorem~\ref{thm:nagata})} \\
  & = e_{R_{\s}}(\s R_{\s}, A_{\s}) & \text{(by (\ref{eqn:rank}))} \\
  & = \sum_{S\cap R=\s} [\kappa(S):\kappa(\s)]e_{A_S}(\s A_S, A_S)\\
\intertext{(where the sum is taken over all primes $S$ of $A$ 
contracting to $\s$ in $R$, by (\ref{eqn:mult2}))}
  & \geq e(\s A_{S_1}, A_{S_1})+e(\s A_{S_2}, A_{S_2}) \\
  & \geq e(S_1 A_{S_1},A_{S_1})+e(S_2 A_{S_2}, A_{S_2}) \\
  & =e(A_{S_1})+e(A_{S_2}).\\
\intertext{By~\cite{nagata:lr} (36.8) the rings $A_{S_1}/PA_{S_1}$ and 
$A_{S_2}/QA_{S_2}$ are analytically unramified.  Therefore, 
Theorem~\ref{thm:nagata} implies that}
e(A_P)+e(A_Q)&>e(A_{S_1})+e(A_{S_2})\geq e(A_P)+e(A_Q).
\end{align*}
This clearly gives a contradiction.
\end{proof}

Given two ideals $I\supseteq J$ in a ring $A$, we say that $J$ is a 
\textit{reduction} of $I$ if there is a positive integer $n$ such that 
$I^{n+1}=JI^n$.  If $(A,\n)$ is local and $I$ is a reduction of $\n$, then 
$e(I,A)=e(A)$.
If $A$ is local with infinite residue field, then 
there is a system of parameters of $A$ that generates a reduction of 
$\n$ (c.f., Northcott and Rees~\cite{northcott:rilr}).
We shall call such a system of parameters a \textit{reductive system of 
parameters} of $A$.
If $A$ has infinite residue field 
then Herrmann, 
Ikeda and Orbanz~\cite{herrmann:ebu} Proposition 10.17 says that 
an element $x$ is part of a reductive system of parameters of $A$ 
if and only if the initial form of $x$ 
in the associated graded ring $\gr_{\n}(A)$ is part of a homogeneous 
system of parameters of $\gr_{\n}(A)$.

The following theorem shows that (ID-2) holds for rings containing a field 
and for a certain class of mixed-characteristic rings.

\begin{thm} \label{thm:id-1}
Let $(A,\n)$ be a quasi-unmixed local ring
and assume that one of the following conditions holds.
\begin{enumerate}
\item $\red{A}$ contains a field, or
\item $A$ is a ring of mixed-characteristic such that the residual 
characteristic $p$ is 
part of a reductive system of parameters of $A$.
\end{enumerate}
Let $P$ and $Q$ be prime ideals in $A$ such 
that both $A/P$ and $A/Q$ are analytically unramified, 
$\sqrt{P+Q}=\n$, and $e(A)<e(A_P)+e(A_Q)$.  Then 
$\dim(A/P)+\dim(A/Q)\leq\dim(A)$.
\end{thm}

In the proof of this result we demonstrate that (ID-2) need 
only be verified for complete domains with infinite residue fields.

\begin{proof}
Step 1.  Pass to the ring $A(X)=A[X]_{\n A[X]}$ to assume that the residue 
field of $A$ is infinite.  Let $\n(X)=\n A(X)$, $P(X)=PA(X)$ and $Q(X)=QA(X)$.  
Then $A(X)$ is a local ring with maximal ideal $\n(X)$.  
If $\red{A}$ contains a field, then the same is true of $\red{A(X)}$.
If the sequence 
$p,x_2,\ldots,x_d$ is a reductive system of parameters of $A$, 
then the same sequence is a reductive system of parameters of $A(X)$ 
since for some $t$, 
\[ \n(X)^{t+1}=\n^{t+1}A(X)=(p,x_2,\ldots,x_d)\n^tA(X).\]
Multiplicities are preserved by property 5 at the beginning of this section.  That is, 
$e(A(X))=e(A)$, $e(A(X)_{P(X)})=e(A_P)$ and $e(A(X)_{Q(X)})=e(A_Q)$.
Both rings $A(X)/P(X)$ 
and $A(X)/Q(X)$ are analytically unramified by Nagata~\cite{nagata:lr} 
(36.8).  Thus, if the result holds in $A(X)$, then
\begin{align*}
\dim(A/P)+\dim(A/Q)&=\dim(A(X)/P(X))+\dim(A(X)/Q(X))\\
&\leq\dim(A(X))=\dim(A) 
\end{align*}
as desired.  

Step 2.  Pass to the completion $(A^*,\n^*)$ to assume that $A$ is complete 
and equidimensional with infinite residue field.  
Let $P^*$ be a prime ideal of $A^*$ that is minimal over $PA^*$ such 
that $\Ht(P^*)=\Ht(P)$, and similarly for $Q^*$.  Since $A/P$ is 
analytically unramified, $PA^*_{P^*}=P^*A^*_{P^*}$.  The fact that the 
extension $A_P\rightarrow A^*_{P^*}$ is flat therefore implies that 
$e(A^*_{P^*})=e(A_P)$ by property 5 at the beginning of this section.  
Similarly, $e(A^*_{Q^*})=e(A_Q)$.  
If $\red{A}$ contains a field, then the same is true of $\red{(A^*)}$.
If the sequence 
$p,x_2,\ldots,x_d$ is a reductive system of parameters of $A$, 
then the same sequence is a reductive system of parameters of $A^*$,
as in the previous step.  If the result 
holds for $A^*$, then it holds for $A$, as in the previous step.

Step 3.  Pass to 
the quotient $A/I$ for a suitably chosen minimal prime $I$, 
to assume that $A$ is a complete domain with 
infinite residue field, and that either $A$ itself contains a field or 
$A$ satisfies hypothesis 2 in the statement of the theorem.  To make this 
reduction, it suffices to verify that 
there is a minimal prime $I$ of $A$ contained in 
$P\cap Q$ such that $e(A/I)<e(A_P/IA_P)+e(A_Q/IA_Q)$.  

First, we show how this gives the desired reduction.  Let $A'=A/I$ with 
prime ideals $\n'=\n A'$, $P'=PA'$ and $Q'=\q A'$.  If $\red{A}$ contains a 
field, then $A'$ contains a field since $A'$ is a quotient of 
$\red{A}$.  If the sequence 
$p,x_2,\ldots,x_d$ is a reductive system of parameters of $A$,
then the same sequence is a reductive system of parameters of $A'$, 
as in Step 2.  Since $A$ is 
equidimensional, $\dim(A')=\dim(A)$.  And since $I\subseteq P\cap Q$, it 
follows that $\dim(A'/P')=\dim(A/P)$ and $\dim(A'/Q')=\dim(A/Q)$.  
Therefore, as in Step 1, we may pass to $A'$.

Now we prove that such a minimal prime $I$ exists.  
Let $\{I_1,\ldots,I_g\}=\min(A)$.  
Suppose that $e(A/I_j)\geq e(A_P/I_jA_P)+e(A_Q/I_jA_Q)$ for every $j$ such 
that $I_j\subseteq P\cap Q$.  (This supposition includes 
the hypothetical possibility that no $I_j$ is contained in $P\cap Q$.)
By Theorem~\ref{thm:nagata}, $e(A_P/I_jA_P)\leq e(A/I_j)$ for every $I_j$ 
contained in $P$, and similarly for $I_j$ contained in $Q$.
The Associativity Formula then implies that
\begin{align*}
e(A_P)+e(A_Q)
  &=\sum_{I_j\subseteq P}e(A_P/I_jA_P)\len(A_{I_j})
    +\sum_{I_j\subseteq Q}e(A_Q/I_jA_Q)\len(A_{I_j}) \\
  &=\sum_{\substack{I_j\subseteq P \\ I_j\not\subseteq Q}}
          e(A_P/I_jA_P)\len(A_{I_j})  
   +\sum_{\substack{I_j\subseteq Q \\ I_j\not\subseteq P}}
          e(A_Q/I_jA_Q)\len(A_{I_j})  \\
  &\quad +\sum_{I_j\subseteq P\cap Q} 
         [e(A_P/I_jA_P)+e(A_Q/I_jA_Q)]\len(A_{I_j}) \\
  &\leq\sum_{\substack{I_j\subseteq P \\ I_j\not\subseteq Q}}
          e(A/I_j)\len(A_{I_j})  
   +\sum_{\substack{I_j\subseteq Q \\ I_j\not\subseteq P}}
          e(A/I_j)\len(A_{I_j})  \\
  &\quad +\sum_{I_j\subseteq P\cap Q} 
         e(A/I_j)\len(A_{I_j}) \\
  &\leq\sum_j e(A/I_j)\len(A_{I_j}) \\
  &=e(A).
\end{align*}
This clearly contradicts the 
assumption that $e(A)<e(A_P)+e(A_Q)$.  Thus, 
there is a minimal prime $I$ of $A$ such that 
$I\subseteq P\cap Q$ and $e(A/I)<e(A_P/IA_P)+e(A_Q/IA_Q)$, as claimed.

Step 4.  We prove the result assuming that $A$ is a complete domain with infinite 
residue field.  
If $A$ contains a field, let $x_1,\ldots,x_d$ be a reductive
system of parameters of $A$, $k$ a 
coefficient field of $A$, $R=k[\![x_1,\ldots,x_d]\!]$, and 
$\m=(x_1,\ldots,x_d)R$.
If $A$ has mixed characteristic, then by assumption there is a reductive system of 
parameters $p,x_2,\ldots,x_d$ of $A$.  By the Cohen Structure Theorem
there is a complete discrete valuation ring $(V,pV)$ contained in $A$ such 
that the induced map on residue fields is an isomorphism.  In this case, let 
$R=V[\![x_2,\ldots,x_d]\!]$ and $\m=(p,x_2,\ldots,x_d)R$.
Note that, in either case $R$ is a 
regular local ring contained in $A$ such that the extension $R\rightarrow A$ is 
module-finite.  
The maximal ideal of $R$ is $\m$ and 
the extension ideal $\m A$ is a reduction of $\n$.
Furthermore,
the induced maps on residue fields is an isomorphism.  
Therefore,
\begin{align*}
e_R(\m,A)&=[A/\n:R/\m]e_A(\m A,A) & \text{(by (\ref{eqn:mult}))} \\
  &=e_A(\n,A)=e(A)<e(A_P)+e(A_Q). 
\end{align*}
If $\p=P\cap R$ and $\q=Q\cap R$, then Lemma~\ref{lem:newrank} implies 
that $\sqrt{\p+\q}=\m$, and Serre's Intersection Theorem implies that
\[ \dim(A/P)+\dim(A/Q)=\dim(R/\p)+\dim(R/\q)\leq\dim(R)=\dim(A).\]
This completes the proof.
\end{proof}

The following theorem is one of our main results and our main application 
of Theorem~\ref{thm:id-1}.

\begin{thm} \label{thm:powers}
Let $(R,\m)$ be a regular local ring containing a field and 
$\p$ and $\q$ prime ideals of $R$ such that $\sqrt{\p+\q}=\m$ and 
$\dim(R/\p)+\dim(R/\q)=\dim(R)$.  Then 
$\p^{(m)}\cap\q^{(n)}\subseteq\m^{m+n}$ for all $m,n\geq 1$.
\end{thm}

\begin{proof}
Fix a nonzero element $f\in \p^{(m)}\cap\q^{(n)}$ and without loss of 
generality, assume that $f\not\in\p^{(m+1)}$ and $f\not\in\q^{(n+1)}$.  
Let $A=R/fR$ with maximal ideal $\n=\m A$ and prime ideals $P=\p A$ and 
$Q=\q A$.  If $f\not\in\m^{m+n}$, then $e(A)<m+n=e(A_P)+e(A_Q)$.  
Theorem~\ref{thm:id-1} implies that
\[ \dim(R/\p)+\dim(R/\q)=\dim(A/P)+\dim(A/Q)\leq\dim(A)=\dim(R)-1 \]
which gives a contradiction.
\end{proof}

It is now straightforward to answer Question~\ref{question:affine} in the 
affirmative, and we do so explicitly in the following corollary.

\begin{cor} \label{cor:affine}
Let $k$ be an algebraically closed field of characteristic 0, and $Y$ and $Z$  
closed subvarieties of $\A^d_k$ that intersect at 
finitely many points, including the origin, and such that 
$\dim(Y)+\dim(Z)=d$.  If $f$ is a nonzero regular 
function on $\A^d_k$ vanishing along $Y$ and $Z$ to orders $m$ and $n$, 
respectively, then $f$ vanishes to order at least $m+n$ at the origin.
\end{cor}

\begin{proof}
Let $\p$ and $\q$ be the prime ideals of $R=k[X_1,\ldots,X_d]$ such that 
$Y=V(\p)$ and $Z=V(\q)$, and let $\m=(X_1,\ldots,X_d)R$.
The aforementioned result of Zariski implies that $f$ is an element of 
$\p^{(m)}\cap\q^{(n)}$.  Passing the the regular local ring $R_{\m}$ we 
have 
\[ f\in\p^{(m)} R_{\m}\cap\q^{(n)}R_{\m}
  =(\p R_{\m})^{(m)}\cap(\q R_{\m})^{(n)}\subseteq(\m R_{\m})^{m+n}
  =\m^{m+n} R_{\m} \]
which implies that $f$ is in the ideal $\m^{(m+n)}=\m^{m+n}$, as desired.
\end{proof}

In the following corollary we state an immediate generalization of 
Theorem~\ref{thm:powers}.  It is easy to see that the corresponding result 
for mixed characteristic will hold once we verify  (SP-2) for such rings.
Recall that a (not 
necessarily local) Noetherian ring is said to be \textit{regular} if the localization 
at every maximal ideal is a regular local ring.

\begin{cor} \label{cor:moreideals}
Let $R$ be a regular ring containing a field with prime ideals 
$\p_1,\ldots,\p_t$.  Let $\m$ be a minimal prime ideal of 
$\p_1+\cdots+\p_t$ and assume that $\Ht(\m)=\sum_i\Ht(\p_i)$.  Then, 
\[ \p_1^{(n_1)}\cap\cdots\cap\p_t^{(n_t)}\subseteq \m^{(\sum_in_i)} \]
for all $n_1,\ldots,n_t\geq 1$.
\end{cor}

One version of Serre's Intersection Theorem states that, since $\m$ is 
minimal over $\sum_i\p_i$, $\Ht(\m)\leq\sum_i\Ht(\p_i)$.  Thus, the 
assumption that $\Ht(\m)=\sum_i\Ht(\p_i)$ simply assures that $\Ht(\m)$ is 
as large as possible.

\begin{proof}
By passing to the localization $R_{\m}$, we assume without loss of 
generality that $R$ is regular local with maximal ideal $\m$.  We prove the result 
in this case by induction on $t$.  When $t=1$, the result is trivial.  
When $t=2$, this is exactly Theorem~\ref{thm:powers}.  Assume now that 
$t\geq 3$.  Fix a prime ideal $\s$ minimal over $\sum_{i\geq 2}\p_i$.  
Then by 
Serre's Intersection Theorem, $\Ht(\s)\leq\sum_{i\geq 2}\Ht(\p_1)$.  
Since $\p_1+\s$ is $\m$-primary, 
$\sum_{i\geq 1}\Ht(\p_1)=\Ht(\m)\leq\Ht(\p_1)+\Ht(\s)$ and it follows that 
$\Ht(\s)=\sum_{i\geq 2}\Ht(\p_1)$.  Therefore, by induction 
\[ \p_1^{(n_1)}\cap\cdots\cap\p_t^{(n_t)}\subseteq
  \p_1^{(n_1)}\cap\s^{(\sum_{i\geq 2}n_i)}\subseteq\m^{(\sum_{i\geq 1}n_i)} \]
as desired.
\end{proof}

In the following corollary, similar to Theorem~\ref{thm:powers}, we give a 
partial verification of (SP-2) in the unramified case of 
mixed-characteristic.

\begin{cor} \label{cor:powers-mixed}
Let $(R,\m,k)$ be an unramified regular local ring of 
mixed-characteristic $p$.  Let 
$\p$ and $\q$ be prime ideals of 
$R$ such that $\sqrt{\p+\q}=\m$ and $\dim(R/\p)+\dim(R/q)=d$ and fix a 
nonzero element $f\in\p^{(m)}\cap\q^{(n)}$.  In the associated graded 
ring $\gr_{\m}(R)$, let the initial forms of $f$ and $p$ be 
denoted by $F$ and $P$, respectively.  
If $F$ is not divisible by $P$, then $f\in\m^{m+n}$.
\end{cor}

\begin{proof}
By passing to the completion $R^*$of the ring $R(X)=R[X]_{\m R[X]}$, we may 
assume without loss of generality that $R$ is complete with infinite residue 
field.  The  
only nontrivial observation is that, if $\gr_{\m}(R)=k[P,X_1,\ldots,X_d]$, 
then the associated graded rings of both $R(X)$ and $R^*$ are 
$k(X)[P,X_1,\ldots,X_d]$.  If $P$ does not divide $F$ in 
$k[P,X_1,\ldots,X_d]$, then the same is true in $k(X)[P,X_1,\ldots,X_d]$.

The remainder of the proof is similar to that of Theorem~\ref{thm:powers}.  The only 
difference is that, in the application of Theorem~\ref{thm:id-1} we must 
verify that $p$ is part of a reductive system of parameters of $A=R/fR$.  
Note that the 
associated graded ring of $A$ is $\gr_{\n}(A)=k[P,X_2,\ldots,X_d]/(F)$.  
The fact that $k$ is infinite and $P$ does not divide $F$ implies that $P$ is part of a 
homogeneous system of parameters of degree 1 for $\gr_{\n}(A)$.  As
noted in the discussion preceding Theorem~\ref{thm:id-1}, this implies that 
$p$ is part of a reductive system of parameters of $A$.  
Thus, Theorem~\ref{thm:id-1} applies.
\end{proof}

\section{Special Cases} \label{sec:easyresults}

We quickly demonstrate three special cases where (SP-2) holds.
First, we prove the case when
one of the ideals is generated by part of a regular system of parameters 
of $R$.

\begin{prop} \label{prop:others}
Let $(R,\m)$ be a regular local ring with prime ideals $\p$ and 
$\q$ such that 
$\sqrt{\p+\q}=\m$ and $\dim(R/\p)+\dim(R/\q)=\dim(R)$.  
If $R/\p$ is regular, 
then $\p^{(m)}\cap \q^{(n)}\subseteq\p^{(m)}\m^n\subseteq\m^{m+n}$ for all $m,n\geq 1$.
\end{prop}

\begin{proof}
Since $R/\p$ is regular, $\p$ 
is generated by elements $x_1,\ldots,x_h$ that form part of a regular 
system of parameters of $R$.  Complete this to a regular system of 
parameters $\x=x_1,\ldots,x_h,\ldots,x_d$ for $R$.  The fact that $\x$ is a regular 
sequence and $R$ is Cohen-Macaulay implies that $\p^{(m)}=\p^m$ for all 
$m\geq 1$.  

It suffices to show that $\p^{(m)}\cap \q^{(n)}\subseteq\m^{m+n}$ for 
all $m,n\geq 1$.  To see that this is sufficient, observe that, because 
$\p$ is generated by part of the regular sequence that generates $\m$, 
\[ \p^{(m)}\m^{n-1}\cap\m^{m+n}=\p^{(m)}\m^n \]
even when $n=1$.  The desired sufficiency now follows by induction on $n$.

By passing to the ring $R(X)=R[X]_{\m R[X]}$ 
as in the proof of Theorem~\ref{thm:id-1}, we may assume without loss of 
generality that the residue field of $R$ is infinite.  We prove that 
$\p^{(m)}\cap\q^{(n)}\subseteq\m^{m+n}$
by induction on $h$.  
If $h=0$ or $h=1$ then $\p$ is principal.  Since $\q$ is 
prime,  it is straightforward to show 
that, in either of these cases, $\p^{(m)}\cap \q^{(n)}=\p^m\q^{(n)}\subseteq\m^{m+n}$.  

Now assume that $h\geq 2$.  Since $R/\m$ is infinite, we 
choose an infinite sequence $a_1,a_2,\ldots\in R$ such that the residues 
$\ol{a_i}$ 
of the $a_i$ modulo $\m$ are all distinct.  Let $y_i=x_1+a_ix_2\in \p$ and 
for each $i$ choose a prime ideal $\q_i$ containing $\q$ and $y_i$ such 
that $\Ht(\q_i)=\Ht(\q)+1$.  The quotient $R_i=R/y_iR$ with prime ideals 
$\p R_i$ and $\q_iR_i$ satisfies the induction hypothesis.  Suppose that $f$ is 
a nonzero element of 
$\p^{(m)}\cap \q^{(n)}\subseteq \p^{(m)}\cap \q_i^{(n)}$.  
Then the residue of $f$ in $R_i$ is in 
$(\p R_i)^{(m)}\cap(\q_iR_i)^{(n)}\subseteq(\m R_i)^{m+n}$ so that 
$f\in y_iR+\m^{m+n}$ for every $i$.  The associated graded ring of $R$ is
the polynomial ring
$\gr_{\m}(r)=k[X_1,\ldots,X_d]$ where $X_i$ is the initial form of $x_i$.
If $f\not\in\m^{m+n}$, then the fact that $f$ is in the ideal $y_iR+\m^{m+n}$ 
for every $i$ implies that the initial form $F$ of $f$ in $\gr_{\m}(R)$ is 
divisible by the initial form of $y_i$.  That is, $F$ is divisible by 
$Y_i=X_1+\ol{a_i}X_2$ for every $i$.  Since the $\ol{a_i}$ are all distinct 
in $R/\m$, the $Y_i$ all generate distinct prime ideals in $k[X_1,\ldots,X_d]$.  
Thus, $F$ has an infinite number of distinct, irreducible factors, 
contradicting the fact that $k[X_1,\ldots,X_d]$ is a unique factorization 
domain.
\end{proof}

In the following proposition, we verify (SP-2) when $\p$ and $\q$ are both 
generated by regular sequences.
We say that an ideal $I$ in a ring $A$ is \textit{unmixed} if 
$\dim(A/P)=\dim(A/I)$ for every 
$P\in\ass(A/I)$.  In this case, 
we define the \textit{$n$th symbolic power} of $I$ to be
\[ I^{(n)}=\cap_{P}(I^nA_{P}\cap A)\]
where the intersection is taken over all associated primes $P$ of $A/I$.  If $I$ 
is prime, this agrees with the standard definition of symbolic powers.  
If $A$ is Cohen-Macaulay and $I$ is generated by a regular sequence, then 
$I$ is unmixed, and for all $n\geq 1$, 
$I^{(n)}=I^n$ and $A/I^n$ is Cohen-Macaulay.

\begin{prop}
Let $(R,\m)$ be a regular local ring with ideals $I$ and $J$ such 
that $\sqrt{I+J}=\m$ and $\dim(R/I)+\dim(R/J)=\dim(R)$.  
If both $I$ and $J$ are generated by regular sequences, then
$I^{(m)}\cap J^{(n)}=I^mJ^n\subseteq \m^{m+n}$ for all $m,n\geq 
1$.
\end{prop}

\begin{proof}
By assumption, $R/I^{(m)}$ and $R/J^{(n)}$ are Cohen-Macaulay 
modules such that $\dim(R/I^{(m)})+\dim(R/J^{(n)})=\dim(R)$.  
By Serre~\cite{serre:alm} Corollary to Theorem 4 in Chapter V, 
we know that
$\tor_i(R/I^{(m)},R/J^{(n)})=0$.  By Rotman~\cite{rotman:iha} 
Corollary 11.27 (iii), 
\[ (I^{(m)}\cap J^{(n)})/I^mJ^n=(I^m\cap J^n)/I^mJ^n=\tor_1(R/I^m,R/J^n)=0 \]
and it follows that $(I^{(m)}\cap J^{(n)})=I^mJ^n\subseteq\m^{m+n}$.
\end{proof}

In the following proposition, we verify (SP-2) for graded primes in a 
polynomial algebra.

\begin{prop} 
Let $p$ be a prime number, $(T,pT)$ a discrete valuation ring,
$R=T[X_1,\ldots,X_d]$ with the standard grading,
and $\m=(p,X_1,\ldots,X_d)R$.  Let $\p$ and $\q$ be 
homogeneous prime ideals of $R$ such that $\sqrt{\p+\q}=\m$ and 
$\dim(R/\p)+\dim(R/\q)=d+1$.  Then 
$\p^{(m)}\cap \q^{(n)}\subseteq\m^{m+n}$ for all $m,n\geq 1$.
\end{prop}

\begin{proof}
By passing to the ring $R(Y)=T(Y)[X_1,\ldots,X_d]$ as in our previous 
arguments, we may assume without 
loss of generality that $T$ 
has infinite residue field.
The ordinary power $\p^m$ is a homogeneous ideal, so that its $\p$-primary 
component $\p^{(m)}$ 
is also homogeneous.  Similarly, $\q^{(n)}$ is also homogeneous.  Therefore, to check the 
desired inclusion, we need only check homogeneous elements.  Let $f$ be a 
homogeneous element of $\p^{(m)}\cap \q^{(n)}$.  Write $f$ as a product 
$f=p^sf_1\cdots f_t$ where each $f_i$ is a nonconstant, irreducible, homogeneous 
polynomial of degree $d_i$.

Because $\p$ and $\q$ are homogeneous and $\sqrt{\p+\q}=\m$,
it is straightforward to show that $p$ is an element of either $\p$ or 
$\q$.  We assume without loss of generality that $p\in\p$.  
Then $p^s\in\p^{(s)}\smallsetminus\p^{(s+1)}$.  By Serre's 
Intersection Theorem applied to the ring $R/pR$, our assumptions imply that $p$ can 
not be in both $\p$ and $\q$.  

Let 
$e_i=\max\{e:f_i\in\p^{(e)}\}$ and $e_i'=\max\{e:f_i\in\q^{(e)}\}$ 
for $i=1,\ldots,t$.
Then, $s+\sum_i e_i\geq m$, $\sum_ie_i'\geq n$, and
each $f_i\in\p^{(e_i)}\cap\q^{(e_i')}$.  If each 
$f_i\in\m^{e_i+e_i'}$ (this is automatic if either $e_i=0$ or $e_i'=0$) then
\[ f=p^sf_1\cdots f_t\in\m^{s+\sum_i(e_i+e_i')}\subseteq\m^{m+n} \]
as desired.  Thus, we may assume without loss of generality that $f$ is 
irreducible and not divisible by $p$.  In particular, $f$ has positive 
degree and some monomial term of $f$ has unit coefficient.  This implies 
that, in the associated graded ring of the localization $R_{\m}$, the 
initial form of $f$ has a monomial term with unit coefficient.  
Therefore, 
by Corollary~\ref{cor:powers-mixed} $f\in\m^{(m+n)}=\m^{m+n}$.
\end{proof}

Department of Mathematics\\
\indent University of Illinois\\
\indent 273 Altgeld Hall,\\
\indent 1409 West Green Street\\
\indent Urbana IL, 61801\\
\indent {\tt ssather@math.uiuc.edu}


\end{document}